\documentclass{amsart}
\usepackage{latexsym,amssymb,amsmath,amsthm,amscd,graphicx,esint}

\makeatletter
\@namedef{subjclassname@2010}{%
  \textup{2010} Mathematics Subject Classification}
\makeatother

\numberwithin{equation}{section}
\newtheorem{theorem}{Theorem}[section]
\newtheorem{lemma}[theorem]{Lemma}

\newtheorem{proposition}[theorem]{Proposition}

\theoremstyle{definition}

\theoremstyle{remark}

\newcommand{\cD}{{\mathcal D}}

\newcommand{\cF}{{\mathcal F}}

\newcommand{\cW}{{\mathcal W}}

\newcommand{\R}{{\mathbb R}}

\newcommand{\Z}{{\mathbb Z}}

\def\al{\alpha}

\def\sg{\sigma}

\def\0{\emptyset}
\def\6{\partial}
\def\8{\infty}

\def\l{\left}
\def\r{\right}
\def\ds{\displaystyle}

\begin{document}

\title[The $n$ linear embedding theorem]{The $n$ linear embedding theorem}
\author[H.~Tanaka]{Hitoshi Tanaka}
\address{Graduate School of Mathematical Sciences, The University of Tokyo, Tokyo, 153-8914, Japan}
\email{htanaka@ms.u-tokyo.ac.jp}
\thanks{
The author is supported by 
the FMSP program at Graduate School of Mathematical Sciences, the University of Tokyo, 
and Grant-in-Aid for Scientific Research (C) (No.~23540187), 
the Japan Society for the Promotion of Science. 
}
\date{}

\begin{abstract}
Let $\sg_i$, $i=1,\ldots,n$, denote 
positive Borel measures on $\R^d$, 
let $\cD$ denote the usual collection of dyadic cubes in $\R^d$ 
and let $K:\,\cD\to[0,\8)$ be a~map. 
In this paper we give a~characterization of the $n$ linear embedding theorem. 
That is, we give a~characterization of the inequality 
$$
\sum_{Q\in\cD}
K(Q)\prod_{i=1}^n\l|\int_{Q}f_i\,d\sg_i\r|
\le C
\prod_{i=1}^n
\|f_i\|_{L^{p_i}(d\sg_i)}
$$
in terms of multilinear Sawyer's checking condition and 
discrete multinonlinear Wolff's potential, 
when $1<p_i<\8$. 
\end{abstract}

\maketitle

\section{Introduction}\label{sec1}
The purpose of this paper is to investigate the $n$ linear embedding theorem. 
We first fix some notations. 
We will denote by $\cD$ the family of all dyadic cubes 
$Q=2^{-k}(m+[0,1)^d)$, 
$k\in\Z,\,m\in\Z^d$. 
Let $K:\,\cD\to[0,\8)$ be a~map and 
let $\sg_i$, $i=1,\ldots,n$, be positive Borel measures on $\R^d$. 
In this paper we give a~necessary and sufficient condition 
for which the inequality 
\begin{equation}\label{1.1}
\sum_{Q\in\cD}
K(Q)\prod_{i=1}^n\l|\int_{Q}f_i\,d\sg_i\r|
\le C
\prod_{i=1}^n
\|f_i\|_{L^{p_i}(d\sg_i)},
\end{equation}
to hold when $1<p_i<\8$. 

For the bilinear embedding theorem, 
in the case $\frac1{p_1}+\frac1{p_2}\ge 1$, 
Sergei Treil gives a~simple proof of the following. 

\begin{proposition}[{\rm\cite[Theorem 2.1]{Tr}}]\label{prp1.1}
Let $K:\,\cD\to[0,\8)$ be a~map and 
let $\sg_i$, $i=1,2$, be positive Borel measures on $\R^d$. 
Let $1<p_i<\8$ and 
$\frac1{p_1}+\frac1{p_2}\ge 1$.
The following statements are equivalent:

\begin{itemize}
\item[{\rm(a)}] 
The following bilinear embedding theorem holds:
$$
\sum_{Q\in\cD}
K(Q)\prod_{i=1}^2\l|\int_{Q}f_i\,d\sg_i\r|
\le c_1
\prod_{i=1}^2
\|f_i\|_{L^{p_i}(d\sg_i)}
<\8;
$$
\item[{\rm(b)}] 
For all $Q\in\cD$, 
$$
\begin{cases}\ds
\l(\int_{Q}\l(\sum_{Q'\subset Q}K(Q')\sg_1(Q')1_{Q'}\r)^{p_2'}\,d\sg_2\r)^{1/p_2'}
\le c_2
\sg_1(Q)^{1/p_1}
<\8,
\\ \ds
\l(\int_{Q}\l(\sum_{Q'\subset Q}K(Q')\sg_2(Q')1_{Q'}\r)^{p_1'}\,d\sg_1\r)^{1/p_1'}
\le c_2
\sg_2(Q)^{1/p_2}
<\8.
\end{cases}
$$
\end{itemize}
\noindent
Moreover,
the least possible $c_1$ and $c_2$ are equivalent.
\end{proposition}
Here, for each $1<p<\8$, 
$p'$ denote the dual exponent of $p$, 
i.e., $p'=\frac{p}{p-1}$, and 
$1_{E}$ stands for the characteristic function of the set $E$. 

Proposition \ref{prp1.1} was first proved 
for $p_1=p_2=2$ in \cite{NaTrVo} 
by the Bellman function method. 
Later in \cite{LaSaUr}, 
this was proved in full generality. 
The checking condition in Proposition \ref{prp1.1} is called 
\lq\lq the Sawyer type checking condition", 
since this was first introduced by Eric~T. Sawyer 
in \cite{Sa1,Sa2}. 

To describe the case 
$\frac1{p_1}+\frac1{p_2}<1$, 
we need discrete Wolff's potential. 

Let $\mu$ and $\nu$ be positive Borel measures on $\R^d$ and 
let $K:\,\cD\to[0,\8)$ be a~map. 
For $p>1$, 
the discrete Wolff's potential 
$\cW^p_{K,\mu}[\nu](x)$ 
of the measure $\nu$ 
is defined by 
$$
\cW^p_{K,\mu}[\nu](x)
:=
\sum_{Q\in\cD}
K(Q)\mu(Q)
\l(
\frac1{\mu(Q)}\sum_{Q'\subset Q}
K(Q')\mu(Q')\nu(Q')
\r)^{p-1}
1_{Q}(x),
\quad x\in\R^d.
$$
The author prove the following. 

\begin{proposition}[{\rm\cite[Theorem 1.3]{Ta1}}]\label{prp1.2}
Let $K:\,\cD\to[0,\8)$ be a~map and 
let $\sg_i$, $i=1,2$, be positive Borel measures on $\R^d$.
Let $1<p_i<\8$ and 
$\frac1{p_1}+\frac1{p_2}<1$.
The following statements are equivalent:

\begin{itemize}
\item[{\rm(a)}] 
The following bilinear embedding theorem holds:
$$
\sum_{Q\in\cD}
K(Q)\prod_{i=1}^2\l|\int_{Q}f_i\,d\sg_i\r|
\le c_1
\prod_{i=1}^2
\|f_i\|_{L^{p_i}(d\sg_i)}
<\8;
$$
\item[{\rm(b)}] 
For $\frac1r+\frac1{p_1}+\frac1{p_2}=1$,
$$
\begin{cases}\ds
\|\cW^{p_2'}_{K,\sg_2}[\sg_1]^{1/p_2'}\|_{L^r(d\sg_1)}
\le c_2<\8,
\\[2mm]\ds
\|\cW^{p_1'}_{K,\sg_1}[\sg_2]^{1/p_1'}\|_{L^r(d\sg_2)}
\le c_2<\8.
\end{cases}
$$
\end{itemize}
\noindent
Moreover,
the least possible $c_1$ and $c_2$ are equivalent.
\end{proposition}

In his excerent survey of the $A_2$ theorem \cite{Hy}, 
Tuomas~P. Hyt\"{o}nen introduces 
another proof of Proposition \ref{prp1.1}, 
which uses the \lq\lq parallel corona" decomposition. 
In this paper, 
following Hyt\"{o}nen's arguments in \cite{Hy}, 
we shall establish the following theorems 
(Theorems \ref{thm1.3} and \ref{thm1.4}). 

\begin{theorem}\label{thm1.3}
Let $K:\,\cD\to[0,\8)$ be a~map and 
let $\sg_i$, $i=1,\ldots,n$, be positive Borel measures on $\R^d$. 
Let $1<p_i<\8$ and 
$\sum_{i=1}^n\frac1{p_i}\ge 1$.
The following statements are equivalent:

\begin{itemize}
\item[{\rm(a)}] 
The following $n$ linear embedding theorem holds:
$$
\sum_{Q\in\cD}
K(Q)\prod_{i=1}^n\l|\int_{Q}f_i\,d\sg_i\r|
\le c_1
\prod_{i=1}^n
\|f_i\|_{L^{p_i}(d\sg_i)}
<\8;
$$
\item[{\rm(b)}] 
For all $j=1,\ldots,n$ and for all $Q\in\cD$, 
$$
\sum_{Q'\subset Q}
K(Q')\sg_j(Q')
\prod_{\substack{i=1 \\ i\ne j}}^n
\l|\int_{Q'}f_i\,d\sg_i\r|
\le c_2
\sg_j(Q)^{1/p_j}
\prod_{\substack{i=1 \\ i\ne j}}^n
\|f_i\|_{L^{p_i}(d\sg_i)}
<\8.
$$
\end{itemize}
\noindent
Moreover,
the least possible $c_1$ and $c_2$ are equivalent.
\end{theorem}

Let the symmetric group $S_n$ be the set of all permutations of 
the set $\{1,\ldots,n\}$, 
that is, 
the set of all bijections 
from the set $\{1,\ldots,n\}$ to itself. 
Let $K:\,\cD\to[0,\8)$ be a~map and 
let $\sg_i$, $i=1,\ldots,n$, be positive Borel measures on $\R^d$. 
Let $1<p_i<\8$ and 
$\sum_{i=1}^n\frac1{p_i}<1$.

Let $\phi\in S_n$. Set 
\begin{align*}
&\frac1{r^{\phi}_1}+\frac1{p_{\phi(1)}}=1,
\\[2mm]
&\frac1{r^{\phi}_2}+\frac1{p_{\phi(1)}}+\frac1{p_{\phi(2)}}=1,
\\ &\qquad\vdots \\[2mm]
&\frac1{r^{\phi}_{n-1}}+\sum_{i=1}^{n-1}\frac1{p_{\phi(i)}}=1,
\\[2mm]
&\frac1r+\sum_{i=1}^n\frac1{p_{\phi(i)}}=1.
\end{align*}
Let, for $Q\in\cD$, 
$$
K^{\phi}_1(Q)
:=
K(Q)\sg_{\phi(1)}(Q)
\l(
\frac1{\sg_{\phi(1)}(Q)}
\sum_{Q'\subset Q}
K(Q')
\prod_{i=1}^n\sg_{\phi(i)}(Q')
\r)^{r^{\phi}_1-1},
$$
let 
$$
K^{\phi}_2(Q)
:=
K^{\phi}_1(Q)\sg_{\phi(2)}(Q)
\l(
\frac1{\sg_{\phi(2)}(Q)}
\sum_{Q'\subset Q}
K^{\phi}_1(Q')
\prod_{i=2}^n\sg_{\phi(i)}(Q')
\r)^{r^{\phi}_2/r^{\phi}_1-1}
$$
and, inductively, 
for $j=3,\ldots,n-1$, let 
$$
K^{\phi}_j(Q)
:=
K^{\phi}_{j-1}(Q)\sg_{\phi(j)}(Q)
\l(
\frac1{\sg_{\phi(j)}(Q)}
\sum_{Q'\subset Q}
K^{\phi}_{j-1}(Q')
\prod_{i=j}^n\sg_{\phi(i)}(Q')
\r)^{r^{\phi}_j/r^{\phi}_{j-1}-1}.
$$

\begin{theorem}\label{thm1.4}
With the notation above, 
the following statements are equivalent:

\begin{itemize}
\item[{\rm(a)}] 
The following $n$ linear embedding theorem holds:
$$
\sum_{Q\in\cD}
K(Q)\prod_{i=1}^n\l|\int_{Q}f_i\,d\sg_i\r|
\le c_1
\prod_{i=1}^n
\|f_i\|_{L^{p_i}(d\sg_i)}
<\8;
$$
\item[{\rm(b)}] 
For all $\phi\in S_n$, 
$$
\l\|\l(
\sum_{Q\in\cD}
K^{\phi}_{n-1}(Q)
1_{Q}
\r)^{1/r^{\phi}_{n-1}}
\r\|_{L^r(d\sg_{\phi(n)})}
\le c_2<\8.
$$
\end{itemize}
\noindent
Moreover,
the least possible $c_1$ and $c_2$ are equivalent.
\end{theorem}

Even though Theorems \ref{thm1.3} and \ref{thm1.4} 
both characterize the same $n$ linear embedding theorem, 
it seems that the characterizations are very different. 
In very recent paper \cite{HaHyLi}, 
Timo~S. H\"{a}nninen, 
Tuomas~P. Hyt\"{o}nen and 
Kangwei Li give a~unified approach saying 
\lq\lq sequential testing" characterization, 
when $n=2,3$. 
Especially, 
our Theorem \ref{thm1.4} with $n=3$ 
is obtained in \cite[Theorem 1.16]{HaHyLi}. 
(An alternative form of another unified characterization has been simultaneously obtained by Vuorinen \cite{Vu}.)
In \cite{Ta2}, the author gives 
a~characterization of the trilinear embedding theorem 
interms of Theorem \ref{thm1.3} and Propositions \ref{prp1.1} and \ref{prp1.2}.

The letter $C$ will be used for constants
that may change from one occurrence to another.

\section{Proof of the necessity}\label{sec2}
In what follows we shall prove the necessity of theorems. 
The necessity of Theorem \ref{thm1.3}, 
that is, 
(b) follows from (a) 
at once if we substitute the test function $f_j=1_{Q}$. 
So, we shall verify the necessity of Theorem \ref{thm1.4}.
We need a~lemma (cf. Lemma 2.1 in \cite{Ta1}).

\begin{lemma}\label{lem2.1}
Let $\sg$ be a~positive Borel measure on $\R^d$. 
Let $1<s<\8$ and 
$\{\al_{Q}\}_{Q\in\cD}\subset[0,\8)$.
Define, for $Q_0\in\cD$, 
\begin{align*}
A_1
&:=
\int_{Q_0}
\l(\sum_{Q\subset Q_0}\frac{\al_{Q}}{\sg(Q)}1_{Q}\r)^s
\,d\sg,
\\
A_2
&:=
\sum_{Q\subset Q_0}
\al_{Q}\l(\frac1{\sg(Q)}\sum_{Q'\subset Q}\al_{Q'}\r)^{s-1},
\\
A_3
&:=
\int_{Q_0}
\sup_{Q\subset Q_0}
\l(\frac{1_{Q}(x)}{\sg(Q)}\sum_{Q'\subset Q}\al_{Q'}\r)^s
\,d\sg(x).
\end{align*}
Then 
$$
A_1\le c(s)A_2,\quad
A_2\le c(s)^{\frac1{s-1}}A_3
\quad\text{and}\quad
A_3\le (s')^sA_1.
$$
Here, 
$$
c(s)
:=
\begin{cases}\ds
s,\quad 1<s\le 2, 
\\ \ds
\l(s(s-1)\cdots(s-k)\r)^{\frac{s-1}{s-k-1}},
\quad 2<s<\8,
\end{cases}
$$
where $k=\lceil s-2 \rceil$ is 
the smallest integer greater than $s-2$. 
\end{lemma}

We will use 
$\fint_{Q}f\,d\sg$ 
to denote the integral average 
$\sg(Q)^{-1}\int_{Q}f\,d\sg$. 
The dyadic maximal operator 
$M_{\cD}^{\sg}$ 
is defined by
$$
M_{\cD}^{\sg}f(x)
:=
\sup_{Q\in\cD}
\frac{1_{Q}(x)}{\sg(Q)}\int_{Q}|f(y)|\,d\sg(y).
$$

Suppose that (a) of Theorem \ref{thm1.4}. 
Then, for $\phi\in S_n$, 
\begin{equation}\label{2.1}
\sum_{Q\in\cD}
K(Q)
\prod_{i=1}^n
\l|\int_{Q}f_{\phi(i)}\,d\sg_{\phi(i)}\r|
\le c_1
\prod_{i=1}^n
\|f_{\phi(i)}\|_{L^{p_{\phi(i)}}(d\sg_{\phi(i)})}.
\end{equation}
Recall that 
$\frac1{r^{\phi}_1}+\frac1{p_{\phi(1)}}=1$. 
By duality, we see that 
$$
\int_{\R^d}\l(
\sum_{Q\in\cD}
K(Q)
\prod_{i=2}^n
\l|\int_{Q}f_{\phi(i)}\,d\sg_{\phi(i)}\r|
1_{Q}
\r)^{r^{\phi}_1}
\,d\sg_{\phi(1)}
\le c_1^{r^{\phi}_1}
\prod_{i=2}^n
\|f_{\phi(i)}\|_{L^{p_{\phi(i)}}(d\sg_{\phi(i)})}^{r^{\phi}_1},
$$
which implies by Lemma \ref{lem2.1} 
\begin{align*}
\lefteqn{
\sum_{Q\in\cD}
K(Q)\sg_{\phi(1)}(Q)
\prod_{i=2}^n
\l|\int_{Q}f_{\phi(i)}\,d\sg_{\phi(i)}\r|
}\\ &\quad\times\l[
\frac1{\sg_{\phi(1)}(Q)}
\sum_{Q'\subset Q}
K(Q')\sg_{\phi(1)}(Q')
\prod_{i=2}^n
\l|\int_{Q'}f_{\phi(i)}\,d\sg_{\phi(i)}\r|
\r]^{r^{\phi}_1-1}
\\ &\le C c_1^{r^{\phi}_1}
\prod_{i=2}^n
\|f_{\phi(i)}\|_{L^{p_{\phi(i)}}(d\sg_{\phi(i)})}^{r^{\phi}_1}.
\end{align*}
It follows from this inequality that 
\begin{align*}
\lefteqn{
\sum_{Q\in\cD}
K^{\phi}_1(Q)
\prod_{i=2}^n
\l|\int_{Q}g_{\phi(i)}\,d\sg_{\phi(i)}\r|
}\\ &=
\sum_{Q\in\cD}
K(Q)\sg_{\phi(1)}(Q)
\prod_{i=2}^n
\l|\int_{Q}g_{\phi(i)}\,d\sg_{\phi(i)}\r|
\l[
\frac1{\sg_{\phi(1)}(Q)}
\sum_{Q'\subset Q}
K(Q')
\prod_{i=1}^n
\sg_{\phi(i)}(Q')
\r]^{r^{\phi}_1-1}
\\ &=
\sum_{Q\in\cD}
K(Q)\sg_{\phi(1)}(Q)
\prod_{i=2}^n
\sg_{\phi(i)}(Q)
\l|\fint_{Q}g_{\phi(i)}\,d\sg_{\phi(i)}\r|^{1/r^{\phi}_1}
\\ &\quad\times
\l[
\frac1{\sg_{\phi(1)}(Q)}
\sum_{Q'\subset Q}
K(Q')\sg_{\phi(1)}(Q')
\prod_{i=2}^n
\sg_{\phi(i)}(Q')
\l|\fint_{Q}g_{\phi(i)}\,d\sg_{\phi(i)}\r|^{1/r^{\phi}_1}
\r]^{r^{\phi}_1-1}
\\ &\le 
\sum_{Q\in\cD}
K(Q)\sg_{\phi(1)}(Q)
\prod_{i=2}^n
\int_{Q}
\l(M_{\cD}^{\sg_{\phi(i)}}g_{\phi(i)}\r)^{1/r^{\phi}_1}
\,d\sg_{\phi(i)}
\\ &\quad\times
\l[
\frac1{\sg_{\phi(1)}(Q)}
\sum_{Q'\subset Q}
K(Q')\sg_{\phi(1)}(Q')
\prod_{i=2}^n
\int_{Q'}
\l(M_{\cD}^{\sg_{\phi(i)}}g_{\phi(i)}\r)^{1/r^{\phi}_1}
\,d\sg_{\phi(i)}
\r]^{r^{\phi}_1-1}
\\ &\le C c_1^{r^{\phi}_1}
\prod_{i=2}^n
\|M_{\cD}^{\sg_{\phi(i)}}g_{\phi(i)}\|_{L^{p_{\phi(i)}/r^{\phi}_1}(d\sg_{\phi(i)})}
\\ &\le C c_1^{r^{\phi}_1}
\prod_{i=2}^n
\|g_{\phi(i)}\|_{L^{p_{\phi(i)}/r^{\phi}_1}(d\sg_{\phi(i)})},
\end{align*}
where we have used the boundedness of dyadic maximal operators. 
Thus, we obtain 
\begin{equation}\label{2.2}
\sum_{Q\in\cD}
K^{\phi}_1(Q)
\prod_{i=2}^n
\l|\int_{Q}f_{\phi(i)}\,d\sg_{\phi(i)}\r|
\le C c_1^{r^{\phi}_1}
\prod_{i=2}^n
\|f_{\phi(i)}\|_{L^{p_{\phi(i)}/r^{\phi}_1}(d\sg_{\phi(i)})}.
\end{equation}
Notice that 
\begin{equation}\label{2.3}
\begin{cases}\ds
\frac{r^{\phi}_{i-1}}{r^{\phi}_i}
+
\frac{r^{\phi}_{i-1}}{p_{\phi(i)}}
=1,\quad
i=2,\ldots,n-1,
\\[4mm]\ds
\frac{r^{\phi}_{n-1}}{r}
+
\frac{r^{\phi}_{n-1}}{p_{\phi(n)}}
=1.
\end{cases}
\end{equation}
By the same manner as the above but starting from \eqref{2.2}, instead of \eqref{2.1}, 
and using \eqref{2.3} with $i=2$, 
we obtain 
$$
\sum_{Q\in\cD}
K^{\phi}_2(Q)
\prod_{i=3}^n
\l|\int_{Q}f_{\phi(i)}\,d\sg_{\phi(i)}\r|
\le C c_1^{r^{\phi}_2}
\prod_{i=3}^n
\|f_{\phi(i)}\|_{L^{p_{\phi(i)}/r^{\phi}_2}(d\sg_{\phi(i)})}.
$$
By being continued inductively until the $n-1$ step,
we obtain 
$$
\sum_{Q\in\cD}
K^{\phi}_{n-1}(Q)
\l|\int_{Q}f_{\phi(n)}\,d\sg_{\phi(n)}\r|
\le C c_1^{r^{\phi}_{n-1}}
\|f_{\phi(n)}\|_{L^{p_{\phi(n)}/r^{\phi}_{n-1}}(d\sg_{\phi(n)})}.
$$
Notice that the last equation of \eqref{2.3}. 
Then by duality 
$$
\l\|
\sum_{Q\in\cD}
K^{\phi}_{n-1}(Q)
1_{Q}
\r\|_{L^{r/r^{\phi}_{n-1}}(d\sg_{\phi(n)})}
\le C c_1^{r^{\phi}_{n-1}}
$$
and, hence, 
$$
\l\|\l(
\sum_{Q\in\cD}
K^{\phi}_{n-1}(Q)
1_{Q}
\r)^{1/r^{\phi}_{n-1}}
\r\|_{L^r(d\sg_{\phi(n)})}
\le C c_1,
$$
which completes the necessity of Theorem \ref{thm1.4}.

\section{Proof of the sufficiency}\label{sec3}
In what follows we shall prove the sufficiency of theorems. 

Let $Q_0\in\cD$ be taken large enough and be fixed. 
We shall estimate the quantity 
\begin{equation}\label{3.1}
\sum_{Q\subset Q_0}
K(Q)
\prod_{i=1}^n\l(\int_{Q}f_i\,d\sg_i\r),
\end{equation}
where $f_i\in L^{p_i}(d\sg_i)$ 
is nonnegative and is supported in $Q_0$. 
We define the collection of principal cubes 
$\cF_i$ for the pair $(f_i,\sg_i)$, 
$i=1,\ldots,n$. Namely, 
$$
\cF_i:=\bigcup_{k=0}^{\8}\cF_i^k,
$$
where 
$\cF_i^0:=\{Q_0\}$,
$$
\cF_i^{k+1}
:=
\bigcup_{F\in\cF_i^k}ch_{\cF_i}(F)
$$
and $ch_{\cF_i}(F)$ is defined by 
the set of all \lq\lq maximal" dyadic cubes $Q\subset F$ such that 
$$
\fint_{Q}f_i\,d\sg_i
>
2\fint_{F}f_i\,d\sg_i.
$$
Observe that
\begin{align*}
\lefteqn{
\sum_{F'\in ch_{\cF_i}(F)}\sg_i(F')
}\\ &\le
\l(2\fint_{F}f_i\,d\sg_i\r)^{-1}
\sum_{F'\in ch_{\cF_i}(F)}
\int_{F'}f_i\,d\sg_i
\\ &\le
\l(2\fint_{F}f_i\,d\sg_i\r)^{-1}
\int_{F}f_i\,d\sg_i
=
\frac{\sg_i(F)}{2},
\end{align*}
which implies 
\begin{equation}\label{3.2}
\sg_i(E_{\cF_i}(F))
:=
\sg_i\l(F\setminus\bigcup_{F'\in ch_{\cF_i}(F)}F'\r)
\ge
\frac{\sg_i(F)}{2},
\end{equation}
where the sets 
$E_{\cF_i}(F)$, $F\in\cF_i$, 
are pairwise disjoint. 
We further define the stopping parents, 
for $Q\in\cD$, 
$$
\begin{cases}\ds
\pi_{\cF_i}(Q)
:=
\min\{F\supset Q:\,F\in\cF_i\},
\\\ds
\pi(Q)
:=
\l(\pi_{\cF_1}(Q),\ldots,\pi_{\cF_n}(Q)\r).
\end{cases}
$$
Then we can rewrite the series in \eqref{3.1} as follows:
$$
\sum_{Q\subset Q_0}
=
\sum_{(F_1,\ldots,F_n)\in(\cF_1,\ldots,\cF_n)}
\sum_{\substack{
Q: \\ \pi(Q)=(F_1,\ldots,F_n)
}}.
$$
We notice the elementary fact that, 
if $P,R\in\cD$, then 
$P\cap R\in\{P,R,\0\}$. 
This fact implies, 
if $\pi(Q)=(F_1,\ldots,F_n)$, 
then 
$$
Q\subset F_{\phi(1)}\subset\cdots\subset F_{\phi(n)}
\quad\text{for some}\quad
\phi\in S_n.
$$
Thus, for fixed $\phi\in S_n$, 
we shall estimate 
\begin{equation}\label{3.3}
\sum_{\substack{
(F_{\phi(i)})\in(\cF_{\phi(i)}):
\\
F_{\phi(1)}\subset\cdots\subset F_{\phi(n)}
}}
\sum_{\substack{
Q: \\ \pi(Q)=(F_{\phi(i)})
}}
K(Q)\prod_{i=1}^n\l(\int_{Q}f_{\phi(i)}\,d\sg_{\phi(i)}\r).
\end{equation}

\paragraph{{\bf Proof of (a) of Theorem \ref{thm1.3}.}}
It follows that, for fixed 
$F_{\phi(n)}\in\cF_{\phi(n)}$, 
\begin{align*}
\lefteqn{
\sum_{F_{\phi(1)}\subset\cdots\subset F_{\phi(n)}}
\sum_{\substack{
Q: \\ \pi(Q)=(F_{\phi(i)})
}}
K(Q)\prod_{i=1}^n\l(\int_{Q}f_{\phi(i)}\,d\sg_{\phi(i)}\r)
}\\ &\le
\l(2\fint_{F_{\phi(n)}}f_{\phi(n)}\,d\sg_{\phi(n)}\r)
\sum_{F_{\phi(1)}\subset\cdots\subset F_{\phi(n)}}
\sum_{\substack{
Q: \\ \pi(Q)=(F_{\phi(i)})
}}
K(Q)\sg_{\phi(n)}(Q)
\prod_{i=1}^{n-1}\l(\int_{Q}f_{\phi(i)}\,d\sg_{\phi(i)}\r).
\end{align*}

We need two observations. 
Suppose that 
$F_{\phi(1)}\subset\cdots\subset F_{\phi(n)}$ 
and $\pi(Q)=(F_{\phi(i)})$. 
Let $i=1,\ldots,n-1$. If 
$F'\in ch_{\cF_{\phi(n)}}(F_{\phi(n)})$ 
satisfies $F'\subset Q$. Then 
\begin{equation}\label{3.4}
\pi_{\cF_{\phi(n)}}\l(\pi_{\cF_{\phi(i)}}(F')\r)
=
\begin{cases}\ds
F_{\phi(n)},
\quad\text{when}\quad
f'\notin\cF_{\phi(i)},
\\\ds
F',
\quad\text{when}\quad
f'\in\cF_{\phi(i)}.
\end{cases}
\end{equation}
By this observation, we define 
$$
ch_{\cF_{\phi(n)}}^{\phi(i)}(F_{\phi(n)})
:=
\{
F'\in ch_{\cF_{\phi(n)}}(F_{\phi(n)}):\,
\text{ $F'$ satisfies \eqref{3.4}}
\}.
$$
We further observe that, 
when $F'\in ch_{\cF_{\phi(n)}}^{\phi(i)}(F_{\phi(n)})$, 
we can regard $f_{\phi(i)}$ as a~constant on $F'$ in the above integrals, 
that is, 
we can replace $f_{\phi(i)}$ by 
$f_{\phi(i)}^{F_{\phi(n)}}$ 
in the above integrals, where 
$$
f_{\phi(i)}^{F_{\phi(n)}}
:=
f_{\phi(i)}1_{E_{\cF_{\phi(n)}}(F_{\phi(n)})}
+
\sum_{F'\in ch_{\cF_{\phi(n)}}^{\phi(i)}(F_{\phi(n)})}
\l(\fint_{F'}f_{\phi(i)}\,d\sg_{\phi(i)}\r)
1_{F'}.
$$

It follows from (b) of Theorem \ref{thm1.3} that 
\begin{align*}
\lefteqn{
\sum_{F_{\phi(1)}\subset\cdots\subset F_{\phi(n)}}
\sum_{\substack{
Q: \\ \pi(Q)=(F_{\phi(i)})
}}
K(Q)\sg_{\phi(n)})(Q)
\prod_{i=1}^{n-1}
\l(\int_{Q}f_{\phi(i)}^{F_{\phi(n)}}\,d\sg_{\phi(i)}\r)
}\\ &\le c_2
\sg_{\phi(n)}(F_{\phi(n)})^{1/p_{\phi(n)}}
\prod_{i=1}^{n-1}
\|f_{\phi(i)}^{F_{\phi(n)}}\|_{L^{p_{\phi(i)}}(d\sg_{\phi(i)})}.
\end{align*}
Thus, we obtain 
$$
\eqref{3.3}
\le C c_2
\sum_{F_{\phi(n)}\in\cF_{\phi(n)}}
\prod_{i=1}^{n-1}
\|f_{\phi(i)}^{F_{\phi(n)}}\|_{L^{p_{\phi(i)}}(d\sg_{\phi(i)})}
\l(\fint_{F_{\phi(n)}}f_{\phi(n)}\,d\sg_{\phi(n)}\r)
\sg_{\phi(n)}(F_{\phi(n)})^{1/p_{\phi(n)}}.
$$
Since 
$\sum_{i=1}^n\frac1{p_{\phi(i)}}\ge 1$, 
we can select the auxiliary parameters 
$s_{\phi(i)}$, $i=1,\ldots,n-1$, 
that satisfy 
$$
\sum_{i=1}^{n-1}\frac1{s_{\phi(i)}}
+
\frac1{p_{\phi(n)}}
=1
\quad\text{and}\quad
1<p_{\phi(i)}\le s_{\phi(i)}<\8.
$$
It follows from H\"{o}lder's inequality with exponents 
$s_{\phi(1)},\ldots,s_{\phi(n-1)},p_{\phi(n)}$ 
that 
\begin{align*}
\eqref{3.3}
&\le C c_2
\prod_{i=1}^{n-1}
\l[
\sum_{F_{\phi(n)}\in\cF_{\phi(n)}}
\|f_{\phi(i)}^{F_{\phi(n)}}\|_{L^{p_{\phi(i)}}(d\sg_{\phi(i)})}^{s_{\phi(i)}}
\r]^{1/s_{\phi(i)}}
\\ &\quad\times
\l[
\sum_{F_{\phi(n)}\in\cF_{\phi(n)}}
\l(\fint_{F_{\phi(n)}}f_{\phi(n)}\,d\sg_{\phi(n)}\r)^{p_{\phi(n)}}
\sg_{\phi(n)}(F_{\phi(n)})
\r]^{1/p_{\phi(n)}}
\\ &\le C c_2
\prod_{i=1}^{n-1}
\l[
\sum_{F_{\phi(n)}\in\cF_{\phi(n)}}
\|f_{\phi(i)}^{F_{\phi(n)}}\|_{L^{p_{\phi(i)}}(d\sg_{\phi(i)})}^{p_{\phi(i)}}
\r]^{1/p_{\phi(i)}}
\\ &\quad\times
\l[
\sum_{F_{\phi(n)}\in\cF_{\phi(n)}}
\l(\fint_{F_{\phi(n)}}f_{\phi(n)}\,d\sg_{\phi(n)}\r)^{p_{\phi(n)}}
\sg_{\phi(n)}(F_{\phi(n)})
\r]^{1/p_{\phi(n)}}
\\[5mm] &=: C c_2
(I_1)\times\cdots\times(I_n),
\end{align*}
where we have used 
$\|\cdot\|_{l^{p_{\phi(i)}}}\ge\|\cdot\|_{l^{s_{\phi(i)}}}$. 

For $(I_n)$, using 
$\sg_{\phi(n)}(F_{\phi(n)})\le 2\sg_{\phi(n)}(E_{\cF_{\phi(n)}}(F_{\phi(n)}))$ 
(see \eqref{3.2}), 
the fact that 
$$
\fint_{F_{\phi(n)}}f_{\phi(n)}\,d\sg_{\phi(n)}
\le
\inf_{y\in F_{\phi(n)}}
M_{\cD}^{\sg_{\phi(n)}}f_{\phi(n)}(y)
$$
and the disjointness of the sets 
$E_{\cF_{\phi(n)}}(F_{\phi(n)})$, 
we have 
\begin{align*}
(I_n)
&\le C
\l[
\sum_{F_{\phi(n)}\in\cF_{\phi(n)}}
\int_{E_{\cF_{\phi(n)}}(F_{\phi(n)})}
\l(M_{\cD}^{\sg_{\phi(n)}}f_{\phi(n)}\r)^{p_{\phi(n)}}\,d\sg_{\phi(n)}
\r]^{1/p_{\phi(n)}}
\\ &\le C
\l[\int_{Q_0}\l(M_{\cD}^{\sg_{\phi(n)}}f_{\phi(n)}\r)^{p_{\phi(n)}}\,d\sg_{\phi(n)}\r]^{1/p_{\phi(n)}}
\le C
\|f_{\phi(n)}\|_{L^{p_{\phi(n)}}(d\sg_{\phi(n)})}.
\end{align*}
It remains to estimate 
$(I_i)$, 
$i=1,\ldots,n-1$. 
It follows that 
\begin{align*}
(I_i)^{p_{\phi(i)}}
&=
\sum_{F_{\phi(n)}\in\cF_{\phi(n)}}
\int_{E_{\cF_{\phi(n)}}(F_{\phi(n)})}f_{\phi(i)}^{p_{\phi(i)}}\,d\sg_{\phi(i)}
\\ &\quad +
\sum_{F_{\phi(n)}\in\cF_{\phi(n)}}
\sum_{F'\in ch_{\cF_{\phi(n)}}^{\phi(i))}(F_{\phi(n)})}
\l(\fint_{F'}f_{\phi(i)}\,d\sg_{\phi(i)}\r)^{p_{\phi(i)}}
\sg_{\phi(i)}(F').
\end{align*}
By the pairwise disjointness of the sets 
$E_{\cF_{\phi(n)}}(F_{\phi(n)})$, 
it is immediate that
$$
\sum_{F_{\phi(n)}\in\cF_{\phi(n)}}
\int_{E_{\cF_{\phi(n)}}(F_{\phi(n)})}f_{\phi(i)}^{p_{\phi(i)}}\,d\sg_{\phi(i)}
\le
\|f_{\phi(i)}\|_{L^{p_{\phi(i)}}(d\sg_{\phi(i)})}^{p_{\phi(i)}}.
$$
For the remaining double sum, 
there holds by the uniqueness of the parent 
\begin{align*}
\lefteqn{
\sum_{F_{\phi(n)}\in\cF_{\phi(n)}}
\sum_{\substack{
F'\in ch_{\cF_{\phi(n)}}(F_{\phi(n)}):
\\
\text{$F'$ satisfies \eqref{3.4}}
}}
\l(\fint_{F'}f_{\phi(i)}\,d\sg_{\phi(i)}\r)^{p_{\phi(i)}}
\sg_{\phi(i)}(F')
}\\ &\le 2
\sum_{F_{\phi(n)}\in\cF_{\phi(n)}}
\sum_{\substack{
F\in\cF_{\phi(i)}:
\\
\pi_{\cF_{\phi(n)}}(F)=F_{\phi(n)}
}}
\sum_{\substack{
F'\in ch_{\cF_{\phi(n)}}(F_{\phi(n)}): 
\\
\pi_{\cF_{\phi(i)}}(F')=F
}}
\l(\fint_{F'}f_{\phi(i)}\,d\sg_{\phi(i)}\r)^{p_{\phi(i)}}
\sg_{\phi(i)}(F')
\\ &\le 2
\sum_{F\in\cF_{\phi(i)}}
\l(2\fint_{F}f_{\phi(i)}\,d\sg_{\phi(i)}\r)^{p_{\phi(i)}}
\sg_{\phi(i)}(F)
\\ &\le C
\|M_{\cD}^{\sg_{\phi(i)}}f_{\phi(i)}\|_{L^{p_{\phi(i)}}(d\sg_{\phi(i)})}^{p_{\phi(i)}}
\le C
\|f_{\phi(i)}\|_{L^{p_{\phi(i)}}(d\sg_{\phi(i)})}^{p_{\phi(i)}}.
\end{align*}
Altogether, we obtain 
$$
\eqref{3.3}
\le C c_2
\prod_{i=1}^n
\|f_{\phi(i)}\|_{L^{p_{\phi(i)}}(d\sg_{\phi(i)})}.
$$
This yields (a) of Theorem \ref{thm1.3}.

\paragraph{{\bf Proof of (a) of Theorem \ref{thm1.4}.}} 
We shall estimate \eqref{3.3} by use of multinonlinear Wolff's potential. 
We first observe that if 
$F_{\phi(i)}\in\cF_{\phi(i)}$, 
$i=1,\ldots,n$, 
satisfy 
$F_{\phi(1)}\subset\cdots\subset F_{\phi(n)}$ 
and, for some $Q\in\cD$, 
$\pi(Q)=(F_{\phi(i)})$, 
then 
\begin{equation}\label{3.5}
\pi_{\cF_{\phi(j)}}(F_{\phi(i)})
=
F_{\phi(j)}
\quad\text{for all}\quad
1\le i<j\le n.
\end{equation}
Fix 
$F_{\phi(i)}\in\cF_{\phi(i)}$, 
$i=1,\ldots,n$, 
that satisfy \eqref{3.5}. Then 
\begin{align*}
\lefteqn{
\sum_{\substack{
Q: \\ \pi(Q)=(F_{\phi(i)})
}}
K(Q)\prod_{i=1}^n\l(\int_{Q}f_{\phi(i)}\,d\sg_{\phi(i)}\r)
}\\ &\le 
\prod_{i=1}^n
\l(2\fint_{F_{\phi(i)}}f_{\phi(i)}\,d\sg_{\phi(i)}\r)
\sum_{\substack{
Q: \\ \pi(Q)=(F_{\phi(i)})
}}
K(Q)\prod_{i=1}^n\sg_{\phi(i)}(Q).
\end{align*}
Recall that 
\begin{equation}\label{3.6}
\begin{cases}\ds
\frac1{r^{\phi}_j}+\sum_{i=1}^j\frac1{p_{\phi(i)}}=1,
\quad j=1,\ldots,n-1,
\\[4mm]\ds
\frac1r+\sum_{i=1}^n\frac1{p_{\phi(i)}}=1.
\end{cases}
\end{equation}
In the following estimates, 
$\sum_{F_{\phi(1)}}$ runs over all 
$F_{\phi(1)}\in\cF_{\phi(1)}$ 
that satisfy \eqref{3.5} 
for fixed 
$F_{\phi(i)}\in\cF_{\phi(i)}$, 
$i=2,\ldots,n$. 
\begin{align*}
\lefteqn{
\sum_{F_{\phi(1)}}
\l(\fint_{F_{\phi(1)}}f_{\phi(1)}\,d\sg_{\phi(1)}\r)
\sum_{\substack{
Q: \\ \pi(Q)=(F_{\phi(i)})
}}
K(Q)\prod_{i=1}^n\sg_{\phi(i)}(Q)
}\\ &\le 
\sum_{F_{\phi(1)}}
\l(\fint_{F_{\phi(1)}}f_{\phi(1)}\,d\sg_{\phi(1)}\r)
\sum_{Q\subset F_{\phi(1)}}
K(Q)\prod_{i=1}^n\sg_{\phi(i)}(Q)
\\ &=
\sum_{F_{\phi(1)}}
\l(\fint_{F_{\phi(1)}}f_{\phi(1)}\,d\sg_{\phi(1)}\r)
\sg_{\phi(1)}(F_{\phi(1)})^{1/p_{\phi(1)}}
\\ &\quad\times
\l(
\fint_{F_{\phi(1)}}
\l(\sum_{Q\subset F_{\phi(1)}}
K(Q)\prod_{i=2}^n\sg_{\phi(i)}(Q)
1_{Q}
\r)\,d\sg_{\phi(1)}\r)
\sg_{\phi(1)}(F_{\phi(1)})^{1/r^{\phi}_1},
\end{align*}
where we have used \eqref{3.6} with $j=1$. 
By H\"{o}lder's inequality, we have further that 
\begin{align*}
\le& 
\l[
\sum_{F_{\phi(1)}}
\l(\fint_{F_{\phi(1)}}f_{\phi(1)}\,d\sg_{\phi(1)}\r)^{p_{\phi(1)}}
\sg_{\phi(1)}(F_{\phi(1)})
\r]^{1/p_{\phi(1)}}
\\ &\quad\times
\l[
\sum_{F_{\phi(1)}}
\l(
\fint_{F_{\phi(1)}}
\l(
\sum_{Q\subset F_{\phi(1)}}
K(Q)\prod_{i=2}^n\sg_{\phi(i)}(Q)
1_{Q}
\r)\,d\sg_{\phi(1)}
\r)^{r^{\phi}_1}
\sg_{\phi(1)}(F_{\phi(1)})
\r]^{1/r^{\phi}_1}.
\end{align*}
By the same way as the estimate of $(I_n)$, 
we see that the last term is majorized by 
$$
C\l(\int_{F_{\phi(2)}}
\l(\sum_{Q\subset F_{\phi(2)}}
K(Q)\prod_{i=2}^n\sg_{\phi(i)}(Q)
1_{Q}
\r)^{r^{\phi}_1}
\,d\sg_{\phi(1)}
\r)^{1/r^{\phi}_1}.
$$
By Lemma \ref{lem2.1}, 
we have further that 
$$
\le C
\l(\sum_{Q\subset F_{\phi(2)}}
K^{\phi}_1(Q)\prod_{i=2}^n\sg_{\phi(i)}(Q)
\r)^{1/r^{\phi}_1}.
$$
By \eqref{2.3}, we notice that 
\begin{equation}\label{3.7}
\frac1{r^{\phi}_i}
+
\frac1{p_{\phi(i)}}
=
\frac1{r^{\phi}_{i-1}},
\quad i=2,\ldots,n-1.
\end{equation}
In the following estimates, 
$\sum_{F_{\phi(2)}}$ runs over all 
$F_{\phi(2)}\in\cF_{\phi(2)}$ 
that satisfy, for fixed 
$F_{\phi(i)}\in\cF_{\phi(i)}$, 
$i=3,\ldots,n$, 
\begin{equation}\label{3.8}
\pi_{\cF_{\phi(j)}}(F_{\phi(i)})
=
F_{\phi(j)}
\quad\text{for all}\quad
2\le i<j\le n.
\end{equation}
There holds 
\begin{align*}
\lefteqn{
\sum_{F_{\phi(2)}}
\l(\fint_{F_{\phi(2)}}f_{\phi(2)}\,d\sg_{\phi(2)}\r)
\times
\l(\sum_{Q\subset F_{\phi(2)}}
K^{\phi}_1(Q)\prod_{i=2}^n\sg_{\phi(i)}(Q)
\r)^{1/r^{\phi}_1}
}\\ &\quad\times
\l(\sum_{F_{\phi(1)}}
\l(\fint_{F_{\phi(1)}}f_{\phi(1)}\,d\sg_{\phi(1)}\r)^{p_{\phi(1)}}
\sg_{\phi(1)}(F_{\phi(1)})
\r)^{1/p_{\phi(1)}}
\\ &=
\sum_{F_{\phi(2)}}
\l(\fint_{F_{\phi(2)}}f_{\phi(2)}\,d\sg_{\phi(2)}\r)
\sg_{\phi(2)}(F_{\phi(2)})^{1/p_{\phi(2)}}
\\ &\quad\times
\l(\fint_{F_{\phi(2)}}
\l(\sum_{Q\subset F_{\phi(2)}}
K^{\phi}_1(Q)\prod_{i=3}^n\sg_{\phi(i)}(Q)
1_{Q}
\r)\,d\sg_{\phi(2)}
\r)^{1/r^{\phi}_1}
\sg_{\phi(2)}(F_{\phi(2)})^{1/r^{\phi}_2}
\\ &\quad\times
\l(\sum_{F_{\phi(1)}}
\l(\fint_{F_{\phi(1)}}f_{\phi(1)}\,d\sg_{\phi(1)}\r)^{p_{\phi(1)}}
\sg_{\phi(1)}(F_{\phi(1)})
\r)^{1/p_{\phi(1)}},
\end{align*}
where we have used \eqref{3.7} with $i=2$.
Recall that \eqref{3.6} with $j=2$. 
Then H\"{o}lder's inequality gives 
\begin{align*}
\le &
\l[\sum_{F_{\phi(2)}}
\l(\fint_{F_{\phi(2)}}f_{\phi(2)}\,d\sg_{\phi(2)}\r)^{p_{\phi(2)}}
\sg_{\phi(2)}(F_{\phi(2)})
\r]^{1/p_{\phi(2)}}
\\ &\times
\l[
\sum_{F_{\phi(2)}}
\sum_{F_{\phi(1)}}
\l(\fint_{F_{\phi(1)}}f_{\phi(1)}\,d\sg_{\phi(1)}\r)^{p_{\phi(1)}}
\sg_{\phi(1)}(F_{\phi(1)})
\r]^{1/p_{\phi(1)}}
\\ &\times
\l[
\sum_{F_{\phi(2)}}
\l(
\fint_{F_{\phi(2)}}
\l(\sum_{Q\subset F_{\phi(2)}}
K^{\phi}_1(Q)\prod_{i=3}^n\sg_{\phi(i)}(Q)
1_{Q}
\r)\,d\sg_{\phi(2)}
\r)^{r^{\phi}_2/r^{\phi}_1}
\sg_{\phi(2)}(F_{\phi(2)})
\r]^{1/r^{\phi}_2}.
\end{align*}
The last term is majorized by 
$$
C\l(\int_{F_{\phi(3)}}
\l(\sum_{Q\subset F_{\phi(3)}}
K^{\phi}_1(Q)\prod_{i=3}^n\sg_{\phi(i)}(Q)
1_{Q}
\r)^{r^{\phi}_2/r^{\phi}_1}
\,d\sg_{\phi(2)}
\r)^{1/r^{\phi}_2}.
$$
By Lemma \ref{lem2.1}, 
we have further that 
$$
\le C
\l(\sum_{Q\subset F_{\phi(3)}}
K^{\phi}_2(Q)\prod_{i=3}^n\sg_{\phi(i)}(Q)
\r)^{1/r^{\phi}_2}.
$$
By being continued inductively until the $n-1$ step,
we obtain 
\begin{align*}
\eqref{3.3}
&\le C
\l[\sum_{F_{\phi(n)}}
\l(\fint_{F_{\phi(n)}}f_{\phi(n)}\,d\sg_{\phi(n)}\r)^{p_{\phi(n)}}
\sg_{\phi(n)}(F_{\phi(n)})
\r]^{1/p_{\phi(n)}}
\\ &\quad\times
\l[
\sum_{F_{\phi(n)}}
\sum_{F_{\phi(n-1)}}
\l(\fint_{F_{\phi(n-1)}}f_{\phi(n-1)}\,d\sg_{\phi(n-1)}\r)^{p_{\phi(n-1)}}
\sg_{\phi(n-1)}(F_{\phi(n-1)})
\r]^{1/p_{\phi(n-1)}}
\\ &\quad\times\vdots
\\ &\quad\times
\l[
\sum_{F_{\phi(n)}}
\sum_{F_{\phi(n-1)}}
\cdots
\sum_{F_{\phi(1)}}
\l(\fint_{F_{\phi(1)}}f_{\phi(1)}\,d\sg_{\phi(1)}\r)^{p_{\phi(1)}}
\sg_{\phi(1)}(F_{\phi(1)})
\r]^{1/p_{\phi(1)}}
\\ &\quad\times
\l[
\sum_{F_{\phi(n)}}
\l(\fint_{F_{\phi(n)}}
\l(\sum_{Q\subset F_{\phi(n)}}
K^{\phi}_{n-1}(Q)
1_{Q}
\r)
\,d\sg_{\phi(n)}
\r)^{r/r^{\phi}_{n-1}}
\sg_{\phi(n)}(F_{\phi(n)})
\r]^{1/r},
\end{align*}
where 
$\sum_{F_{\phi(n)}}$ 
runs over all 
$F_{\phi(n)}\in\cF_{\phi(n)}$ 
and $\sum_{F_{\phi(k)}}$, 
$k=3,\ldots,n-1$, 
runs over all 
$F_{\phi(k)}\in\cF_{\phi(k)}$ 
that satisfy, for fixed 
$F_{\phi(i)}$, $i=k+1,\ldots,n$, 
\begin{equation}\label{3.9}
\pi_{\cF_{\phi(j)}}(F_{\phi(i)})
=
F_{\phi(j)}
\quad\text{for all}\quad
k\le i<j\le n.
\end{equation}
The last term is majorized by 
$$
C\l(\int_{Q_0}
\l(\sum_{Q\subset Q_0}
K^{\phi}_{n-1}(Q)
1_{Q}
\r)^{r/r^{\phi}_{n-1}}
\,d\sg_{\phi(n)}
\r)^{1/r}
\le c_2.
$$
It follows from \eqref{3.5}, \eqref{3.8}, \eqref{3.9} and 
the uniqueness of the parents that 
\begin{align*}
\lefteqn{
\l[
\sum_{F_{\phi(n)}}
\sum_{F_{\phi(n-1)}}
\cdots
\sum_{F_{\phi(i)}}
\l(\fint_{F_{\phi(i)}}f_{\phi(i)}\,d\sg_{\phi(i)}\r)^{p_{\phi(i)}}
\sg_{\phi(i)}(F_{\phi(i)})
\r]^{1/p_{\phi(i)}}
}\\ &\le 
\l[
\sum_{F_{\phi(n)}}
\sum_{\substack{
F_{\phi(i)}\in\cF_{\phi(i)}:
\\
\pi_{\cF_{\phi(n)}}(F_{\phi(i)})=F_{\phi(n)}
}}
\l(\fint_{F_{\phi(i)}}f_{\phi(i)}\,d\sg_{\phi(i)}\r)^{p_{\phi(i)}}
\sg_{\phi(i)}(F_{\phi(i)})
\r]^{1/p_{\phi(i)}}
\\ &=
\l[
\sum_{F_{\phi(i)}\in\cF_{\phi(i)}}
\l(\fint_{F_{\phi(i)}}f_{\phi(i)}\,d\sg_{\phi(i)}\r)^{p_{\phi(i)}}
\sg_{\phi(i)}(F_{\phi(i)})
\r]^{1/p_{\phi(i)}}
\\[5mm] &\le C
\|f_{\phi(i)}\|_{L^{p_{\phi(i)}}(d\sg_{\phi(i)})}.
\end{align*}
Altogether, we obtain 
$$
\eqref{3.3}
\le C c_2
\prod_{i=1}^n
\|f_{\phi(i)}\|_{L^{p_{\phi(i)}}(d\sg_{\phi(i)})}.
$$
This yields (a) of Theorem \ref{thm1.4}.

\end{document}